\documentclass[12pt]{amsart}
\usepackage{amssymb}
\usepackage{epsfig}  		% For postscript
\usepackage{epic,eepic}       % For epic and eepic output from xfig
\usepackage{showkeys}
\usepackage[left=1cm,right=2cm]{geometry}
\newtheorem{thm}{Theorem}[section]

\theoremstyle{definition}
\newtheorem{definition}[thm]{Definition}
\newtheorem{example}[thm]{Example}

\theoremstyle{remark}
\theoremstyle{definition}

\numberwithin{equation}{section}
\begin{document}

%% The title of the paper goes here.  Edit to your title.
\title{$I$-convergent triple difference sequence spaces using sequence of modulus function}
%% Now edit the following to give your name and address:

\author{Tanweer Jalal}
\address{Department of  Mathematics,  National Institute of Technology, Hazratbal, Srinagar-190006, Jammu and Kashmir, India.}
\email{tjalal@nitsri.net}
%%
%% If there is another author uncomment and edit the following.
%%

\author{Ishfaq Ahmad Malik}
\email{ishfaq$\_$2phd15@nitsri.net}
%\urladdr{www.math.sc.edu/$\sim$second}

%%
%% If there are three of more authors they are added in the obvious
%% way. 
%%

%%%
%%% The following is for the abstract.  The abstract is optional and
%%% if not used just delete, or comment out, the following.
%%%

\begin{abstract}
The main objective of this paper is to  introduce classes of $I$-convergent triple difference sequence spaces, $c_{0I}^{3}(\Delta,\digamma)$, $c_{I}^{3}(\Delta,\digamma)$, $\ell_{\infty I}^{3}(\Delta,\digamma)$, $M_{I}^{3}(\Delta,\digamma)$  and $M_{0I}^{3}(\Delta,\digamma)$, by using sequence of modulii function $\digamma=(f_{pqr})$. We also study some algebraic and topological properties of these new sequence spaces.
\end{abstract}

\maketitle
%%%%%%%%%%%%%%%%%%%%%%%%%%%%%%%%%%%%%%%%%%%%%%%%%%%%%%%%%%%%%%%%%%%%%%
\parindent=0mm\vspace{0.05in}

\parindent=0mm\vspace{0.05in}
\section{Introduction}
%%%%%%%%%%%%%%%%%%%%%%%%%%%%%%%%%%%%%%%%%%%%%%%%%%%%%%%%%%%%%%%%%%%%%%
$~~~~$A triple sequence (real or complex) is a function $x:\mathbb{N}\times \mathbb{N} \times \mathbb{N}\rightarrow \mathbb{R}(\mathbb{C})$, where $\mathbb{N},\mathbb{R} \text{ and } \mathbb{C} $ are the set of natural numbers, real numbers, and complex numbers respectively. We denote by $\omega^{'''}$ the class of all complex triple sequence $(x_{pqr}) $, where $p,q,r \in \mathbb{N}$. Then under the coordinate wise addition and scalar multiplication $\omega^{'''}$ is a linear space. 
A triple sequence can be represented by a matrix, in case of double sequences we write in the form of a square. In case of triple sequence it will be in the form of a box in three dimensions. \\

\parindent=5mm\vspace{0.00in}
The different types of notions of triple sequences and their statistical convergence were introduced and investigated initially by Sahiner et. al [16]. Later Debnath et.al [1,2], Esi et.al [3,4,5], and many others authors have studied it further and obtained various results. Kizmaz [11] introduced the notion of difference sequence spaces, he defined the difference sequence spaces $ \ell_\infty(\Delta)~,~ c(\Delta)$ and $c_0(\Delta)$ as follows.
$$Z(\Delta)=\{x=(x_k)\in \omega \colon (\Delta x_k)\in Z\} $$
for $Z=\ell_\infty\,,\, c \,\text{ and } c_0$\\
Where $\Delta x=(\Delta x_k)=(x_k-x_{k+1})$  and $\Delta^{0}x_k=x_k$ for all $k\in \mathbb{N}$ \\
The difference operator on triple sequence is defined as
\begin{align*}
\Delta x_{mnk}=x_{mnk}-&x_{(m+1)nk}-x_{m(n+1)k}-x_{mn(k+1)}+x_{(m+1)(n+1)k}\\
&+x_{(m+1)n(k+1)} +x_{m(n+1)(k+1)}-x_{(m+1)(n+1)(k+1)}
\end{align*}
and $\Delta_{mnk}^{0}=( x_{mnk} ) $.

\parindent=5mm\vspace{0.00in}
Statistical convergence was introduced by Fast [6] and later on it was studied by Fridy [7-8] from the sequence space point of view and linked it with summability theory. The notion of statistical convergent double sequence was introduced by Mursaleen and Edely [14].

\parindent=5mm\vspace{0.00in}
$I$-convergence is a generalization of the statistical convergence. Kostyrko et. al. [12] introduced the notion of $I$-convergence of real sequence and studied its several properties. Later Jalal [9-10], Salat et.al [15]  and many other researchers contributed in its study. Tripathy and Goswami [19] extended this concept in probabilistic normed space using triple difference sequences of real numbers. Sahiner and Tripathy [17] studied $I$-related properties in triple sequence spaces and showed some interesting results. Tripathy[18] extended the concept in $I$-convergent double sequence and later  Kumar [13] obtained some results on $I$-convergent double sequence.

\parindent=5mm\vspace{0.00in}
In this paper we have defined $I$-convergent triple difference sequence spaces, $c_{0I}^{3}(\Delta,\digamma)$, $c_{I}^{3}(\Delta,\digamma)$, $\ell_{\infty I}^{3}(\Delta,\digamma)$, $M_{I}^{3}(\Delta,\digamma)$  and $M_{0I}^{3}(\Delta,\digamma)$, by using sequence of modulii function $\digamma=(f_{pqr})$ and also studied some algebraic and topological properties of these new sequence spaces. 

\section{Definitions and preliminaries}
\begin{definition}
Let $X\not=\phi$. A class $I\subset 2^X$ (Power set of $X$) is said to be an ideal in $X$ if the following conditions holds good:
\begin{enumerate}
\item[(i)] $I$  is additive that is if $A,B \in I$ then $A\cup B \in I$;
\item[(ii)] $I$  is hereditary that is if $A\in I$, and $B\subset A$ then $B\in I $.	
\end{enumerate}
$I$ is called non-trivial ideal if $X\not \in I$
\end{definition}

\begin{definition}[16] A triple sequence $( x_{pqr})$ is said to be convergent to $L$ in Pringsheim's sense if for every $\epsilon>0$, there exists $\bold{N}\in \mathbb{N}$ such that $$|x_{pqr}-L|<\epsilon~~~~~~~\text{  whenever  } ~~~p\geq \bold{N},q\geq \bold{N},r\geq \bold{N}$$
and write as $\lim_{p,p,r\rightarrow\infty}x_{pqr}=L$.
\end{definition}
\parindent=0mm\vspace{0.00in}
\textbf{Note:} A triple sequence is convergent in Pringsheim's sense may not be bounded [16].\\

\parindent=0mm\vspace{0.00in}
\textbf{Example} Consider the sequence $(x_{pqr})$ defined by
$$x_{pqr}=\left\{\begin{matrix}
p+q& \text{for all}~ p=q~\text{and} ~r=1\\
{1\over p^2qr}& \text{otherwise} 
\end{matrix}\right.$$
Then $x_{pqr} \rightarrow 0$ in Pringsheim's sense but is unbounded.
\begin{definition} A triple sequence $(x_{pqr})$ is said to be $I$-convergence to a number $L$ if for every $\epsilon>0$ , 
$$\left\{(p, q, r) \in \mathbb{N}\times \mathbb{N} \times \mathbb{N} : |x_{pqr} - L|) \geq \epsilon\right\} \in  I .$$
In this case we write $I - \lim x_{pqr} = L$ .
\end{definition}
\begin{definition} A triple sequence $(x_{pqr})$ is said to be $I$-null if $L = 0$. In this case we write $I- \lim x_{pqr} = 0$ .
\end{definition}
\begin{definition} [16] A triple sequence $( x_{pqr})$ is said to be Cauchy sequence if for every $\epsilon>0,$ there exists $\bold{N}\in \mathbb{N}$ such that 
$$|x_{pqr}-x_{lmn}|<\epsilon~~~~\text{  whenever  } ~p\geq l\geq \bold{N},q\geq m\geq \bold{N},r\geq n\geq \bold{N}$$
\end{definition}
\begin{definition} A triple sequence $( x_{pqr})$ is said to be $I-$Cauchy sequence if for every $\epsilon>0,$ there exists $\bold{N}\in \mathbb{N}$ such that 
$$\{(p,q,r)\in \mathbb{N}\times \mathbb{N} \times \mathbb{N}:|x_{pqr}-a_{lmn}|\geq \epsilon\}\in I$$
  $\text{  whenever  } ~~~p\geq l\geq \bold{N},q\geq m\geq \bold{N},r\geq n\geq \bold{N}$
\end{definition}
\begin{definition} [16] A triple sequence $( x_{pqr})$ is said to be bounded if there exists $M>0,$ such that $|x_{pqr}|<M$ for all $p,q,r\in\mathbb{N}$.
\end{definition}
\begin{definition} A triple sequence $( x_{pqr})$ is said to be $I-$bounded if there exists $M>0,$ such that $\{(p, q, r) \in \mathbb{N}\times \mathbb{N} \times \mathbb{N} : |x_{pqr} | \geq M\} \in  I$ for all $p,q,r\in\mathbb{N}$.
\end{definition}
\begin{definition} A triple sequence space $E$ is said to be solid if $(\alpha_{pqr}x_{pqr})\in E$ whenever $(x_{pqr}) \in  E$ and for all sequences $(\alpha_{pqr})$ of scalars with $|\alpha_{pqr}| \leq 1$, for all $p, q, r \in \mathbb{N}$ .
\end{definition}
\begin{definition} Let $E$ be a triple sequence space and $x=(x_{pqr})\in  E$. Define the set $S(x)$ as  
$$S(x)=\left\{\left(x_{\pi(pqr)} \right):\pi\text{ is a permutations of } \mathbb{N} \right\} $$
If $S(x)\subseteq E$ for all $x\in E$, then $E$ is said to be symmetric.
\end{definition}
\begin{definition} A triple sequence space $E$ is said to be convergence free  if $(y_{pqr})\in  E$ whenever $(x_{pqr}) \in  E$ and $x_{pqr} = 0$ implies $y_{pqr} = 0$ for all $p,q,r \in \mathbb{N}$.
\end{definition}

\begin{definition} A triple sequence space $E$ is said to be sequence algebra if $x\cdot y \in E$ , whenever $x=(x_{pqr}) \in E$ and $y=(y_{pqr}) \in E$, that is product of any two sequences is also in the space. 
\end{definition}
\begin{definition} A function $f:[0,\infty)\rightarrow [0,\infty)$ is called a modulus function
if it satisfies the following conditions 
\begin{enumerate}
\item[(i)] $f(x)=0$  if and only if $x=0$.
\item[(ii)] $f(x+y)\leq f(x)+f(y)$ for all $x\geq 0$ and $y\geq 0$.
\item[(iii)] $f$ is increasing.
\item[(iv)] $f$ is continuous from the right at 0.
\end{enumerate}
Since $|f(x)-f(y)|\leq f(|x-y|)$, it follows from condition (4) that $f$ is continuous on $[0,\infty)$. Furthermore, from condition (2) we have $f (nx) \leq nf (x)$, for all $n \in \mathbb{N}$, and so\\
$f(x)=f\big(nx({1\over n})\big)\leq nf\big( {x\over n}\big) $.\\
 Hence ${1\over n}f(x)\leq f({x\over n})$
for  all $n\in \mathbb{N}$
\end{definition}

We now define the following sequence spaces 
\begin{align*}
c_{0I}^{3}(\Delta,\digamma)=&\left\{x\in \omega^{'''}:I-\lim f_{pqr}(|\Delta x_{pqr}|)=0\right\}\\
c_{I}^{3}(\Delta,\digamma)=&\left\{x\in \omega^{'''}:I-\lim f_{pqr}(|\Delta x_{pqr}-b|)=0, ~\text{ for some }~b  \right\} \\
\ell_{\infty I}^{3}(\Delta,\digamma)=&\left \{x\in \omega^{'''}:\sup_{p,q,r\in\mathbb{N}} f_{pqr}(|\Delta x_{pqr}|)=0\right\}\\
 M_{I}^{3}(\Delta,\digamma)=& c_{I}^{3}(\Delta,\digamma) \cap \ell_{\infty I}^{3}(\Delta,\digamma) \\ 
M_{0I}^{3}(\Delta,\digamma)=&c_{0I}^{3}(\Delta,\digamma)\cap \ell_{\infty I}^{3}(\Delta,\digamma)
\end{align*}

\section{Algebraic and Topological Properties of the new Sequence spaces}
\begin{thm} The triple difference sequence spaces $c_{0I}^{3}(\Delta,\digamma), c_{I}^{3}(\Delta,\digamma), \ell_{\infty I}^{3}(\Delta,\digamma), M_{I}^{3}(\Delta,\digamma)$ and  $ M_{0I}^{3}(\Delta,\digamma)$ all are linear for the sequence of modulii $\digamma=( f_{pqr} ) $.  
\end{thm}
\begin{proof} We shall prove it for the sequence space $c_{I}^{3}(\Delta,\digamma)$, for the other spaces, it can be established similarly.\\
Let $x=( x_{pqr} ),y=(y_{pqr} ) \in c_{I}^{3}(\Delta,\digamma)$ and $\alpha, \beta \in\mathbb{R}$ such that $|\alpha|\leq 1$ and $|\beta|\leq 1$, then
\begin{align*}
&I-\lim f_{pqr}(|\Delta x_{pqr}-b_1|)=0,~~\text{for some  } b_1\in \mathbb{C}\\
&I-\lim f_{pqr}(|\Delta y_{pqr}-b_2|)=0,~~\text{for some  } b_2\in \mathbb{C}
\end{align*}
Now for a given $\epsilon>0$ we set 
$$X_1=\left\{ (p,q,r)\in \mathbb{N}\times\mathbb{N}\times \mathbb{N}:f_{pqr}(|\Delta x_{pqr}-b_1|)>{\epsilon\over 2}\right\}\in I \,\qquad \qquad \qquad(2.1)$$
$$X_2=\left\{ (p,q,r)\in \mathbb{N}\times\mathbb{N}\times \mathbb{N}:f_{pqr}(|\Delta y_{pqr}-b_2|)>{\epsilon\over 2}\right\}\in I  \,\qquad \qquad \qquad (2.2)$$
Since $f_{pqr}$ is a modulus function, so it is non-decreasing and convex, hence we get
\begin{align*}
f_{pqr}(|(\alpha\Delta x_{pqr}+\beta\Delta y_{pqr})-&(\alpha b_1+\beta b_2)|)=
f_{pqr}(|(\alpha\Delta x_{pqr}-\alpha b_1)+(\beta \Delta y_{pqr}-\beta b_2)|)\\
&\leq f_{pqr}(|\alpha||\Delta x_{pqr}-b_1|)+f_{pqr}(|\beta||\Delta y_{pqr}-b_2|)\\
&=|\alpha|f_{pqr}(|\Delta x_{pqr}-b_1|)+|\beta|f_{pqr}(|\Delta y_{pqr}-b_2|)\\
&\leq f_{pqr}(|\Delta x_{pqr}-b_1|)+f_{pqr}(|\Delta y_{pqr}-b_2|)
\end{align*}
From (2.1) and (2.2) we can write
$$\left\{ (p,q,r)\in \mathbb{N}\times\mathbb{N}\times \mathbb{N}:f_{pqr}(|(\alpha\Delta x_{pqr}+\beta\Delta y_{pqr})-(\alpha b_1+\beta b_2)|)>\epsilon\right\}\subset X_1\cup X_2 $$
Thus $\alpha x+\beta y \in c_{I}^{3}(\Delta,\digamma)$\\
This completes the proof.
\end{proof}
\begin{thm} The triple difference sequence $x=(x_{pqr}) \in M_{I}^{3}(\Delta,\digamma)$ is $I$-convergent if and only if for every $\epsilon>0	$ there exists $I_\epsilon,J_\epsilon,K_\epsilon\in \mathbb{N}$ such that 
$$\left\{ (p,q,r)\in \mathbb{N}\times\mathbb{N}\times \mathbb{N}:f_{pqr}(|\Delta x_{pqr}-\Delta x_{I_\epsilon J_\epsilon K_\epsilon}|)\leq \epsilon\right\}\in M_{I}^3(\Delta,\digamma) $$
\end{thm}
\begin{proof} Let $b=I-\lim\Delta x.$ Then we have 
$$A_\epsilon=\left\{ (p,q,r)\in \mathbb{N}\times\mathbb{N}\times \mathbb{N}:f_{pqr}(|\Delta x_{pqr}-b|)\leq {\epsilon\over 2}\right\}\in M_{I}^3(\Delta,\digamma) \qquad \text{ for all} , \epsilon>0. $$
Next fix $I_\epsilon,J_\epsilon,K_\epsilon\in A_\epsilon$ then we have 
$$|\Delta x_{pqr}-\Delta x_{I_\epsilon J_\epsilon K_\epsilon}|\leq 
|\Delta x_{pqr}-b|+|b-\Delta x_{I_\epsilon J_\epsilon K_\epsilon}|\leq {\epsilon\over2}+{\epsilon\over 2}=\epsilon\, \text{for all, } \, p,q,r\in A_\epsilon$$
Thus 
$$\left\{ (p,q,r)\in \mathbb{N}\times\mathbb{N}\times \mathbb{N}:f_{pqr}(|\Delta x_{pqr}-\Delta x_{I_\epsilon J_\epsilon K_\epsilon}|)\leq \epsilon\right\}\in M_{I}^3(\Delta,\digamma) $$
Conversely suppose that 
$$\left\{ (p,q,r)\in \mathbb{N}\times\mathbb{N}\times \mathbb{N}:f_{pqr}(|\Delta x_{pqr}-\Delta x_{I_\epsilon J_\epsilon K_\epsilon}|)\leq \epsilon\right\}\in M_{I}^3(\Delta,\digamma) $$
we get  $\left\{ (p,q,r)\in \mathbb{N}\times\mathbb{N}\times \mathbb{N}:f_{pqr}(|\Delta x_{pqr}-\Delta x_{I_\epsilon J_\epsilon K_\epsilon}|)\leq \epsilon\right\}\in M_{I}^3(\Delta,\digamma) $, for all $\epsilon>0$ . \\
Then $\text{given}~\epsilon>0$ we can find the set\\
$$B_\epsilon=\left\{ (p,q,r)\in \mathbb{N}\times\mathbb{N}\times \mathbb{N}:\Delta x_{pqr}\in [\Delta x_{I_\epsilon J_\epsilon K_\epsilon}-\epsilon\,,\, \Delta x_{I_\epsilon J_\epsilon K_\epsilon}+\epsilon]\right\}\in M_{I}^3(\Delta,F)$$
Let $J_\epsilon=[\Delta x_{I_\epsilon J_\epsilon K_\epsilon}-\epsilon\,,\,\Delta x_{I_\epsilon J_\epsilon K_\epsilon}+\epsilon ]$ if $\epsilon>0 $ is fixed than $B_\epsilon \in M_{I}^3(\Delta,F)$ as well as $B_{\epsilon\over 2} \in M_{I}^3(\Delta,\digamma)$. \\
Hence $B_\epsilon\cap  B_{\epsilon\over 2} \in M_{I}^3(\Delta,\digamma)$ \\
Which gives $J=J_\epsilon\cap  J_{\epsilon\over 2}\not=\phi$ that is $\left\{ (p,q,r)\in \mathbb{N}\times\mathbb{N}\times \mathbb{N}:\Delta x_{pqr}\in \mathbb{N}\right\} \in M_{I}^3(\Delta,\digamma)$\\
Which implies $\text{diam}~ J\leq \text{diam} ~J_\epsilon $\\
where the diam of $J$ denotes the the length of interval $J$.\\
Now by the principal of induction a sequence of closed interval can be found 
$$J_\epsilon=I_0 \supseteq  I_1\supseteq I_2 \supseteq\cdots \supseteq I_s \supseteq \cdots \cdots$$
with the help of the property that $\text{diam}~ I_s\leq {1\over 2}\text{diam} ~I_{s-1} $, for $s=1,2,\cdots$ and \\
$\left\{ (p,q,r)\in \mathbb{N}\times\mathbb{N}\times \mathbb{N}:\Delta x_{pqr}\in I_{pqr}\right\}$ for $(p,q,r=1,2,3\cdots)$\\
Then there exists a $\xi \in \cap I_s$ where  $s\in \mathbb{N}$ such that $\xi=I-\lim \Delta x$\\
So that $f_{pqr}(\xi)=I-\lim f_{pqr}(\Delta x)$ therefore $b=I-\lim f_{pqr}(\Delta x)$.\\
Hence the proof is complete.
\end{proof}

\begin{thm}
The $\digamma=(f_{pqr})$ be a sequence of modulus functions then the inclusions $c_{0I}^{3}(\Delta,\digamma)\subset c_{I}^{3}(\Delta,\digamma)\subset \ell_{\infty I}^{3}(\Delta,\digamma) $ holds .
\end{thm}
\begin{proof}
The inclusion $c_{0I}^{3}(\Delta,\digamma)\subset c_{I}^{3}(\Delta,\digamma)$ is obvious.\\
We prove $c_{I}^{3}(\Delta,\digamma)\subset \ell_{\infty I}^{3}(\Delta,\digamma) $, let $x=( x_{pqr})\in c_{I}^{3}(\Delta,\digamma)$ then there exists $b\in \mathbb{C}$ such that 
$I-\lim f_{pqr}(|\Delta x_{pqr}-b|)=0$ ,\\ 
Which gives $f_{pqr}(|\Delta x_{pqr}|)\leq f_{pqr}(|\Delta x_{pqr}-b|)+f_{pqr}(|b|) $\\
On taking supremum over $p,q$ and $r$ on both sides gives \\
$x=(x_{pqr})\in \ell_{\infty I}^{3}(\Delta,\digamma)$\\
Hence the inclusion $c_{0I}^{3}(\Delta,\digamma)\subset c_{I}^{3}(\Delta,\digamma)\subset \ell_{\infty I}^{3}(\Delta,\digamma) $ holds.
\end{proof}
\begin{thm}
The triple difference sequence $c_{0I}^{3}(\Delta,\digamma)$ and $ M_{0I}^{3}(\Delta,\digamma) $ are solid.
\end{thm}	
\begin{proof}
We prove the result for $c_{0I}^{3}(\Delta,\digamma)$. \\
Consider $x=( x_{pqr})\in c_{0I}^{3}(\Delta,\digamma) $, then $I-\lim_{p,q,r} f_{pqr}(|\Delta x_{pqr}|)=0 $\\
Consider a sequence of scalar $(\alpha_{pqr})$ such that $|\alpha_{pqr}|\leq 1$ for all $ p,q,r \in \mathbb{N}$. \\
Then we have 
\begin{align*}
I-\lim_{p,q,r} f_{pqr}(|\Delta \alpha_{pqr}( x_{pqr})|)&\leq I-|\alpha_{pqr}|\lim_{p,q,r} f_{pqr}(|\Delta x_{pqr}|)\\&\leq I-\lim_{p,q,r} f_{pqr}(|\Delta x_{pqr}|)\\ &=0
\end{align*}
Hence $I-\lim_{p,q,r} f_{pqr}(|\Delta \alpha_{pqr} x_{pqr}|)=0 $ for all $p,q,r \in \mathbb{N}$\\
Which gives $(\alpha_{pqr} x_{pqr})\in c_{0I}^{3}(\Delta,\digamma) $\\
Hence the sequence space $c_{0I}^{3}(\Delta,\digamma) $ is solid.\\
The result for $M_{0I}^{3}(\Delta,\digamma)$ can be similarly proved. 
\end{proof}
\begin{thm}
The triple difference sequence spaces $c_{0I}^{3}(\Delta,\digamma)$ , $c_{I}^{3}(\Delta,\digamma)$ , $\ell_{\infty I}^{3}(\Delta,\digamma)$ , $M_{I}^{3}(\Delta,\digamma)$ and $M_{0I}^{3}(\Delta,\digamma)$ are sequence algebras.
\end{thm}
\begin{proof}
We prove the result for $c_{0I}^{3}(\Delta,\digamma)$.\\
Let $x=( x_{pqr}), y=(y_{pqr})\in c_{0I}^{3}(\Delta,\digamma)$\\
Then we have $I-\lim f_{pqr}(|\Delta x_{pqr}|)=0 $  and $I-\lim f_{pqr}(|\Delta y_{pqr}|)=0 $\\
and $I-\lim f_{pqr}(|\Delta (x_{pqr}\cdot y_{pqr})|)=0 $ as
\begin{align*}
\Delta (x_{pqr}\cdot y_{pqr})=x_{pqr}\cdot y_{pqr}-&x_{(p+1)qr}\cdot y_{(p+1)qr}-
x_{p(q+1)r}\cdot y_{p(q+1)r}-x_{pq(r+1)}\cdot y_{pq(r+1)}+x_{(p+1)(q+1)r}\cdot y_{(p+1)(q+1)r}+\\
&x_{(p+1)q(r+1)}\cdot y_{(p+1)q(r+1)}+x_{p(q+1)(r+1)}\cdot y_{p(q+1)(r+1)}-x_{(p+1)(q+1)(r+1)}\cdot y_{(p+1)(q+1)(r+1)}
\end{align*}
It implies that $x\cdot y \in c_{0I}^{3}(\Delta,\digamma)$\\
Hence the proof.\\
The result can be proved for the spaces $c_{I}^{3}(\Delta,\digamma)$ , $\ell_{\infty I}^{3}(\Delta,\digamma)$ , $M_{I}^{3}(\Delta,\digamma)$ and $M_{0I}^{3}(\Delta,\digamma)$ in the same way.
 \end{proof}
\begin{thm}
In general the sequence spaces $c_{0I}^{3}(\Delta,\digamma)$ , $c_{I}^{3}(\Delta,\digamma)$ and $\ell_{\infty I}^{3}(\Delta,\digamma)$ are not convergence free.
\end{thm}
\begin{proof} We prove the result for the sequence space $c_{I}^{3}(\Delta,\digamma)$ using an example.  
\begin{example} Let $I=I_f$ define the triple sequence $x=(x_{pqr})$ as \\ 
 $$ x_{pqr}=\left\{\begin{matrix}0 & \text{ if } p=q=r\\ 1& \text{otherwise}
\end{matrix}\right.$$
Then if $f_{pqr}(x)=x_{pqr}~~\forall ~p,q,r \in\mathbb{N}$, we have  $x=( x_{pqr})\in c_{I}^{3}(\Delta,\digamma)$.\\
Now define the sequence $y=y_{pqr} $ as 
 $$ y_{pqr}=\left\{\begin{matrix}0 & \text{ if } r \text{ is odd , and }~p,q\in\mathbb{N} \\ lmn& \text{otherwise} \end{matrix}\right.$$
Then for $f_{pqr}(x)=x_{pqr}~~\forall ~p,q,r \in\mathbb{N}$, it is clear that  $y=( y_{pqr})\not  \in c_{I}^{3}(\Delta,\digamma)$ \\
Hence the sequence spaces $c_{I}^{3}(\Delta,\digamma)$ is not convergence free.
\end{example}
The space $c_{I}^{3}(\Delta,\digamma)$ and $\ell_{\infty I}^{3}(\Delta,\digamma)$ are not convergence free in general can be proved in the same fashion.
\end{proof}
\begin{thm} In general the triple difference sequences $c_{0I}^{3}(\Delta,\digamma)$  and  $c_{I}^{3}(\Delta,\digamma)$ are not symmetric if $I$ is neither maximal nor $I=I_f$. 
\end{thm}
\begin{proof} We prove the result for the sequence space $c_{0I}^{3}(\Delta,\digamma)$ using an example.  
\begin{example} Define the triple sequence $x=(x_{pqr})$ as \\ 
 $$ x_{pqr}=\left\{\begin{matrix}0 & \text{ if } r=1,~\text{ for all } p,q\in \mathbb{N}\\ 1& \text{otherwise}
\end{matrix}\right.$$
Then if $f_{pqr}(x)=x_{pqr}~~\forall ~p,q,r \in\mathbb{N}$, we have  $x=( x_{pqr})\in c_{0I}^{3}(\Delta,\digamma)$.\\
Now if $x_{\pi(pqr)} $ be a rearrangement of $x=(x_{pqr})$ defined as \\
$$x_{\pi(pqr)}=\left\{\begin{matrix}
1 &\text{for }~p,q,r \text{  even} \in K\\ 0& \text{otherwise} \end{matrix}\right. $$
Then $\{x_{\pi(p,q,r)}\}\not\in c_{0I}^{3}(\Delta,\digamma)$  as $\Delta x_{\pi(pqr)}=1$  \\
Hence the sequence spaces $c_{0I}^{3}(\Delta,\digamma)$ is not symmetric in general.
\end{example}
The space $c_{I}^{3}(\Delta,\digamma)$ is not symmetric in general can be proved in the same fashion.
\end{proof}

\end{document}